\documentclass[a4paper, 
10pt
]{article}

\usepackage[a4paper,left=3cm,right=3cm,top=2.5cm,bottom=2.5cm]{geometry}
\usepackage{subfig}
\usepackage{pgfplots}
\usepackage{pgfplotstable}
\usepackage{mathtools}
\usepackage{multicol}
\usepackage{comment}
\usepackage{booktabs}
\pgfplotsset{compat=1.5}
\usepackage{amssymb}
\usepackage{url}
\usepackage{bm}
\usepackage{siunitx}
\providecommand{\norm}[1]{\lVert#1\rVert}

\usepackage{pdfsync}
\usepackage{float}
\usepackage{tabularx}
\usepackage{enumerate}
\usepackage{array}
\usepackage{xspace}
\usepackage{tikz}
\usepackage{tikzsymbols}
\usetikzlibrary{calc,trees,positioning,arrows,chains,shapes.geometric,%
    decorations.pathreplacing,decorations.pathmorphing,shapes,%
    matrix,shapes.symbols, decorations.markings, patterns,fit}

\usepackage{authblk}

\usepackage{hyperref}
\usepackage[draft,inline,marginclue]{fixme}

\usepackage[utf8]{inputenc}
\usepackage{pgfplots}
\pgfplotsset{compat=newest}
\usepgfplotslibrary{groupplots}
\usepgfplotslibrary{dateplot}

\FXRegisterAuthor{mt}{amt}{MT}

\FXRegisterAuthor{gs}{ags}{GS}

\FXRegisterAuthor{nd}{and}{ND}

\newcommand{\mupar}{\ensuremath{\boldsymbol{\mu}}}
\newcommand{\etapar}{\ensuremath{\boldsymbol{\eta}}}
\newcommand{\domain}{\Omega(\bm \mu)}
\newcommand{\R}{\mathbb{R}}
\newcommand{\sss}{v}
\newcommand{\sv}{\mathbf{v}}

\newcommand{\Am}{\boldsymbol{\mathbf{A}}}
\newcommand{\ATm}{\boldsymbol{\mathbf{\tilde{A}}}}
\newcommand{\Xm}{\boldsymbol{\mathbf{X}}}
\newcommand{\Ym}{\boldsymbol{\mathbf{Y}}}
\newcommand{\Um}{\boldsymbol{\mathbf{U}}}
\newcommand{\Sm}{\boldsymbol{\mathbf{\Sigma}}}
\newcommand{\Vm}{\boldsymbol{\mathbf{V}}}
\newcommand{\PHIm}{\boldsymbol{\mathbf{\Phi}}}
\newcommand{\LAMm}{\boldsymbol{\mathbf{\Lambda}}}
\newcommand{\xv}{\boldsymbol{\mathbf{x}}}

\usepackage{ifpdf}

\ifpdf
\usepackage{pdfsync}
\fi

\begin{document}

\title{Enhancing CFD predictions in shape design problems by model and
parameter space reduction} 

\author[]{Marco~Tezzele\footnote{marco.tezzele@sissa.it}}
\author[]{Nicola~Demo\footnote{nicola.demo@sissa.it}}
\author[]{Giovanni~Stabile\footnote{giovanni.stabile@sissa.it}}
\author[]{Andrea~Mola\footnote{andrea.mola@sissa.it}}
\author[]{Gianluigi~Rozza\footnote{gianluigi.rozza@sissa.it}}

\affil{Mathematics Area, mathLab, SISSA, via Bonomea 265, I-34136
  Trieste, Italy}

\maketitle

\begin{abstract}
In this work we present an advanced computational pipeline for the
approximation and prediction of the lift coefficient of a parametrized
airfoil profile. The non-intrusive reduced order method is based on dynamic
mode decomposition (DMD) and it is coupled with dynamic active
subspaces (DyAS) to enhance
the future state prediction of the target function and reduce the
parameter space dimensionality. The pipeline is based on high-fidelity
simulations carried out by the application of finite volume method
for turbulent flows, and automatic mesh morphing through radial basis
functions interpolation technique. The proposed pipeline is able
  to save $1/3$ of the overall computational resources thanks to the
  application of DMD. Moreover exploiting DyAS and performing the
  regression on a lower dimensional space results in the reduction of
  the relative error in the approximation of the time-varying lift
  coefficient by a factor $2$ with respect to using only the DMD.
\end{abstract}

\tableofcontents

\section{Introduction}
\label{sec:intro}

Reduced order modeling (ROM) is nowadays a quite popular and consolidated technique,
applied to several fields of engineering and computational science thanks to
the remarkable computational gain granted for the solution of the governing equations.
The ROM goal is in fact that of reducing the dimension of the studied system
without altering some important properties of the original problem. This typically
results in more efficient, time saving computations. Among other fields, ROM methods are frequently
and successfully applied to problems governed by parametric partial differential
equations (PDEs), for which many solutions of the same PDE in correspondence with
different parameters are required. This paradigm is for example
encountered in the context of parametric optimal control problems, uncertainty
quantification, and shape optimization.

Model reduction for PDEs has been historically obtained in different ways.
In some cases, very successful reduced models have been obtained at the
level of the governing equations, based on physical
considerations. This is for instance the case of the potential flow theory in the
fluid dynamics field. In other cases, the reduction can be introduced at the
discretization level, as is the case, for instance, for the Boundary Element Method
used in structural analysis, fluid mechanics, electro-magnetism and acoustics studies.
In the case in which parametric PDEs are considered, a possible approach to obtain
efficient reduced order models is to sample the
solution manifold by creating a solutions database corresponding to
different parameters, using a high-dimensional discretization, then
combine the latter to identify the intrinsic lower dimension of the
problem. For parametric reduced order models
see~\cite{hesthaven2016certified,quarteroni2014reduced,morhandbook2019},
while for a more applications oriented overview we
suggest~\cite{tezzele2018ecmi,rozza2018advances,salmoiraghi2016advances}. 

For parametric time-dependent problems, a proper orthogonal decomposition
approach can be applied to reduce the dimensionality of the system, as
in~\cite{GeoStaRoBlu2018,HiStaMoRo2019}. In this work we
propose a novel data-driven approach for parametric dynamical systems, combining 
dynamic mode decomposition (DMD) with active subspaces (AS) property. These
two relatively new methodologies provide a simplification of the
dynamical system, and an analysis of the input parameter space of a
given target function, respectively. Exploiting AS property we are able
to obtain an estimation of the importance of the parameters of such
function, as well as a reduction in the number of parameters. Moreover
the methods are equation-free, being based only on input/output
couples and do not make assumptions on the underlying governing equations.

We define a generic scalar output $\sss(\mupar, t) \in \R$ that depends both on
time $t$ and on the parameters of the model $\mupar \in
\mathbb{D}\subset \mathbb{R}^k$, with $k$ denoting the dimension of the
parameter space. We denote the state of the parametric system 
at time $t$ with $\sss_t(\mupar) \in \R$. The solution manifold in
time is approximated using the DMD in order to obtain an 
approximation of the linear map $A$ defined as:
\begin{equation}
\sss_{t+1}(\mupar) = A(\sss_t(\mupar)).
\label{eq:dmdintro}
\end{equation}
It is easy to note that using~\eqref{eq:dmdintro} we have the possibility to
forecast a generic future state of the parametric system.

To numerically compute the linear operator $A$, we need to sample the parameter
space $\mathbb{D}$, and for each time store the quantity of interest for each
parametric configuration. Formally, considering a set of parameter samples with
dimension $Ns$, the discrete vector referring to the system state at time $t$ results:
\begin{equation}
\sv_t = \begin{bmatrix} \sss_t({\mupar_1}) & \dotsc &
  \sss_t({\mupar_{Ns}})\end{bmatrix}^T \in \R^{Ns}. 
\end{equation}
Collecting several time states $\sv_i(\mupar)$ for $i = 1, \dotsc, m$,
we compute the operator $\mathbf{A}$ with a best-fit approach such that
$\sv_{t+1} \approx \mathbf{A}\sv_t$.  Once computed the future
prevision, we are able to exploit the relation between the input parameters
$\mupar_i$ and the related outputs $\sss_\text{future}(\mupar_i)$ to approximate
the output for any new parameter. In this work we use a Gaussian Process
Regression (GPR)~\cite{williams2006gaussian,guo2018reduced}, but any regression or
interpolation method can be used. We underline that the chosen regression model
has to be fitted for any forecasted time we want to analyse.

The high dimensionality in the parameter space may incur on the
  inability to solve many-query problems with sufficiently high
  fidelity, thus causing a decrease in the accuracy of the solution
  approximation. For this reason we couple the regression with the
AS property in order to perform a sensitivity analysis of function $\sss_t(\mupar)$.
AS indeed is able to provide an approximation $g$ of a scalar
function $f$, where the input parameters of $g$ are a linear combination of the
original parameters of $f$. The coefficients of such combination give
information about the importance of the original parameters. In this work, we
use this information to reduce the dimension of the parameter space --- in
which we build the regression --- by not considering the parameters whose AS
coefficients are smaller than a certain threshold, that is they are almost zero.

The developed methodology is tested on an aeronautics application given by the
flow past an airfoil profile. As output of interest we considered the lift coefficient
and the parameters vector $\mupar$ describes geometrical transformations according
to the morphing technique proposed in~\cite{hicks1978wing}. The fluid
dynamics problem is described using the incompressible Navier--Stokes equations
with turbulence modeling. These are discretized using a finite volume
approximation. The deformed meshes corresponding to different input
parameters are automatically obtained exploiting a Radial Basis
Function (RBF) mesh morphing technique. 

This work is structured as follows: in~\autoref{sec:problem} we present
the general parametric problem over which we apply the proposed numerical
pipeline, providing some information about the geometrical
deformation. In~\autoref{sec:dmd} and~\autoref{sec:as} we 
present the DMD and AS methods, respectively, while in~\autoref{sec:results} we show
the numerical setting of the problem and the results obtained. Finally
in~\autoref{sec:the_end} we propose some final remarks and highlight
possible future developments.

\section{The parametric problem}
\label{sec:problem}
Let be given the unsteady incompressible Navier-Stokes equations
described in an Eulerian framework on a parametrized space-time domain
$Q(\bm \mu) = \domain \times [0,T] \subset
\mathbb{R}^d\times\mathbb{R}^+, \textrm{ }d=2,3$ with the vectorial
velocity field $\bm{u}:Q(\bm \mu) \to \mathbb{R}^d$, and the scalar
pressure field $p:Q(\bm \mu) \to \mathbb{R}$ such that:
\begin{equation}
\label{eq:navstokes}
\begin{cases}
\bm{u_t}+ \bm{\nabla} \cdot (\bm{u} \otimes \bm{u})- \bm{\nabla} \cdot
2 \nu \bm{\nabla^s} \bm{u}=-\bm{\nabla}p &\mbox{ in } Q(\bm \mu),\\
\bm{\nabla} \cdot \bm{u}=\bm{0} &\mbox{ in } Q(\bm \mu),\\
\bm{u} (t,x) = \bm{f}(\bm{x}) &\mbox{ on } \Gamma_{\text{in}} \times [0,T],\\
\bm{u} (t,x) = \bm{0} &\mbox{ on } \Gamma_{0}(\bm \mu) \times [0,T],\\ 
(\nu\nabla \bm{u} - p\bm{I})\bm{n} = \bm{0} &\mbox{ on } \Gamma_{\text{out}} \times [0,T],\\ 
\bm{u}(0,\bm{x})=\bm{k}(\bm{x}) &\mbox{ in } Q(\bm \mu)_0,\\            
\end{cases}
\end{equation}
holds.  Here, $\Gamma = \Gamma_{\text{in}} \cup \Gamma_{0} \cup \Gamma_{\text{out}}$
is the boundary of $\domain$ and it is composed by three different parts
$\Gamma_{\text{in}}$, $\Gamma_{\text{out}}$ and $\Gamma_0(\bm \mu)$ that indicate,
respectively, inlet boundary, outlet boundary, and physical walls. The
term $\bm{f}(\bm{x})$ depicts the stationary non-homogeneous boundary
condition, whereas $\bm{k}(\bm{x})$ denotes the initial condition for
the velocity at $t=0$. Shape changes are applied to the domain $\Omega$,
and in particular to its boundary $\Gamma_0(\bm \mu)$ corresponding to the airfoil
wall. Such shape modifications are associated to numerical parameters contained in the
vector $\bm \mu \in \mathbb{R}^k$ which, in the numerical
examples shown in this work has dimension $k=10$. As said, the only portion of 
the domain boundary subject to shape parametrization is the physical wall 
of the airfoil $\Gamma_0(\bm \mu)$, which in the undeformed configuration corresponds to
the 4-digits, NACA~4412 wing profile~\cite{abbott2012theory,jacobs1933characteristics}. To alter such
geometry, we adopt the shape parametrization and morphing technique proposed
in~\cite{hicks1978wing}, where $k$ shape functions are added to the
airfoil profiles. Let $y_u$, and $y_l$ be the upper and lower ordinates
of a NACA profile, respectively. We express the deformation of such
coordinates as
\begin{align}
  y_u &= \overline{y_u} + \sum_{i=1}^{5} c_i r_i ,\\
  y_l &= \overline{y_l} - \sum_{i=1}^{5} d_i r_i ,
\end{align}
where the bar denotes the reference undeformed state, which is
the NACA~4412 profile.

\begin{figure}[h!]
\centering
\includegraphics[width=0.55\textwidth]{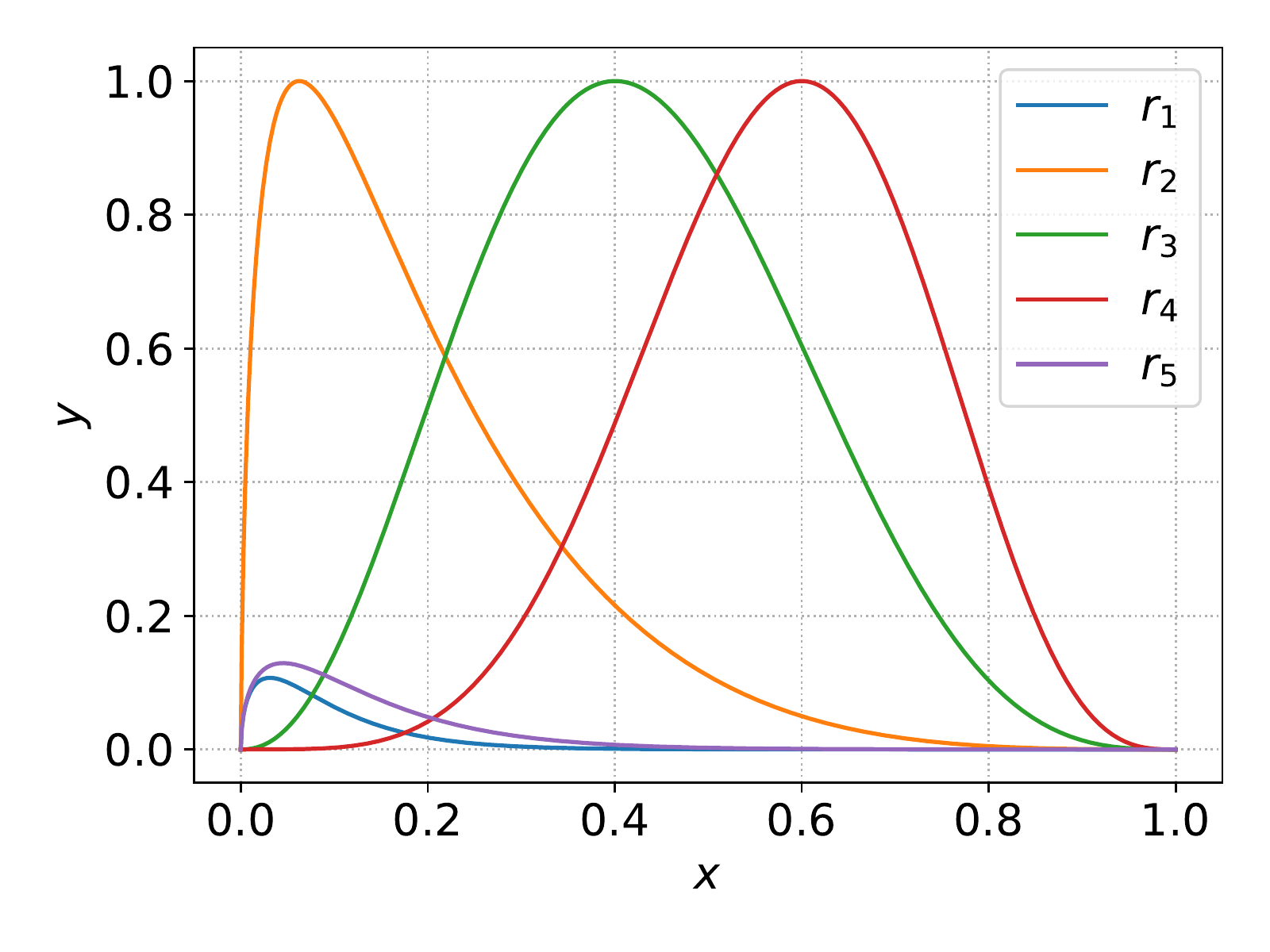}
\caption{Airfoil shape functions with respect to the profile
  abscissa. The leading edge corresponds to $x=0$.}
\label{fig:bump_functions}
\end{figure}

The parameters $\mupar \in \mathbb{D} \subset \mathbb{R}^{10}$ are
the weights coefficients, $c_i$ and $d_i$,  associated with the
shape functions $r_i$. The range of each parameter will be specified
in~\autoref{sec:results}. The explicit formulation of the shape functions
can be found in~\cite{hicks1978wing}, we report them in~\autoref{fig:bump_functions}.

After the reference profile is deformed, we also apply the same
morphing to the mesh coordinates by using a radial basis functions
(RBF) interpolation
method~\cite{buhmann2003radial,morris2008cfd,manzoni2012model}. With
this approach the movement $\bm s$ of all the points which do not belong to
the moving boundaries is approximated by an interpolatory radial basis
function:
\begin{equation}
\label{eq:rbfFunction}
\bm s (\bm x) = \sum_{i=1}^{N_b} \beta_i \xi (||\bm x - \bm x_{b_i}||)+q(\bm x),
\end{equation}
where $\bm{x}_{b_i}$ are the coordinates of points for which we know
the boundary displacements, for this particular case the points
located on the wing surface. $N_b$ is the number of control points on
the boundary, $\xi$ is a given basis function, $q(\bm{x})$ is a
polynomial. The coefficients $\beta_i$ and the polynomial $q(\bm{x})$
are obtained by the imposition of interpolation conditions
\begin{equation}
\bm{s}(\bm{x}_{b_i}) = \bm{d}_{b_i},
\end{equation}
where $\bm{d}_{b_i}$ is the displacement value at the boundary points
and by the additional requirement: 
\begin{equation}
\sum_{i=1}^{N_b}\beta_i q(\bm{x}_{b_i}) = 0.
\end{equation}
In the present case, we select basis functions for which it is
possible to use linear polynomials $q(\bm{x})$. For more informations
concerning the selection of the order of
polynomials see~\cite{Beckert2001}. Finally the values of the coefficients
$\beta_i$ and the coefficients $\delta_i$ of the linear polynomials
$q$ can be obtained by solving the linear problem: 
\begin{equation}
\label{eq:RBF_system}
\begin{bmatrix}
\bm{d}_b \\ 0
\end{bmatrix}
=
\begin{bmatrix}
\bm{M}_{b,b} & P_b \\ P_b^T & 0  
\end{bmatrix}
\begin{bmatrix}
\bm{\beta} \\ \bm{\delta}
\end{bmatrix},
\end{equation}
where $\bm{M}_{b,b} \in \mathbb{R}^{N_b\times N_b}$ is a matrix
containing the evaluation of the basis functions $\xi_{b_i b_j} = \xi
(\norm{\bm{x}_{b_i}-\bm{x}_{b_j}})$, and $\bm{P}_b \in
\mathbb{R}^{N_b\times (d+1)}$ is a matrix where $d$ is the spatial
dimension. Each row of this matrix, that contains the coordinates of
the boundary points, is given by $\mbox{row}_i(\bm P_b)
= \begin{bmatrix}1 & \bm{x_{b_i}}\end{bmatrix}$. Once the system of
\eqref{eq:RBF_system} is solved one can obtain the displacement of
all the internal points using the RBF interpolation:
\begin{equation}
\bm d_{\text{in}_i} = \bm s(\bm x_{\text{in}_i}),
\end{equation}
where $\bm x_{\text{in}_i}$ are the coordinates of the internal grid
points. The computation of the displacement of the grid points entails
the resolution of a dense system of equations that has dimension $N_b
+ d + 1$. Usually, the number of boundary points $N_b$ is much smaller
than the number of grid points $N_h$.

\section{Dynamical systems approximation by dynamic mode decomposition}
\label{sec:dmd}
Dynamic mode decomposition (DMD) is an emerging reduced order method proposed
by Schmid in~\cite{schmid2010dynamic} for the analysis of dynamical systems.
Approximating the linear infinite-dimensional Koopman
operator~\cite{koopman1931hamiltonian}, DMD decomposes the original system into
few main features, the so called DMD modes, that evolve linearly in time, even
if the original system has nonlinear behaviour. This means that, other than
individuating recurrent patterns in the evolution of the system, DMD provides a
real-time midcast/forecast of the output of interest.
An important advantage of such method is the complete data-driven nature: the
algorithm relies only on the system output, without the necessity of any
information regarding the model or equations used.

Dynamic mode decomposition has been successfully employed in naval
hull shape optimization pipelines~\cite{demo2018shape}, for online
real-time acquisitions in a wind tunnel
experiment~\cite{zhang2019online}, and in
meteorology~\cite{bistrian2015improved}, among others. We also mention
the higher order DMD extension~\cite{le2017higher,levega2017higher}. 

In the following paragraph, we provide just an algorithmic overview
of the method. For an exhaustive explanation of DMD, its
applicability, and possible extensions, we
suggest~\cite{kutz2016dynamic,brunton2019data}. 

We define the linear operator $\Am$ such that
\begin{equation}
\label{eq:operator}
\xv_{k+1} = \Am \xv_{k},
\end{equation}
where $\xv_{k+1} \in \R^N$ and $\xv_{k} \in \R^N$ are the vectors containing
the system outputs at two sequential instants. Thus, the operator $\Am: \R^N \to
\R^N$ expresses the dynamics of the system. In order to construct it using only
data, we need to collect $m \leq N+1$ equispaced in time outputs $\xv_i$ for $i = 1,
\dotsc, m$ --- from now on called {\it snapshots} --- then arrange them in two matrices:
$\mathbf{X} = \begin{bmatrix} \mathbf{x}_1 & \dotsc &
\mathbf{x}_{m-1}\end{bmatrix}$ and $\mathbf{Y} = \begin{bmatrix} \mathbf{x}_2 &
\dotsc & \mathbf{x}_{m} \end{bmatrix}$. Since the corresponding columns in
$\Xm$ and $\Ym$ are sequential snapshots, we are able to
use~\eqref{eq:operator} to represent the relationship between $\mathbf{X}$
and $\mathbf{Y}$, such that $\Ym = \Am\Xm$. We can find such operator
  by using the relation $\Am = \Ym \Xm^\dagger$, where $^\dagger$ refers to the
  Moore-Penrose pseudo-inverse. We exploit the singular value
  decomposition to compute such pseudo-inverse, due to its
  computational efficiency and accuracy, as in the following:
% Minimizing the error $\|\Ym - \Am\Xm\|$
% we obtain the linear operator, which still has very large dimension,
% especially when the studied system requires a fine discretization. To reduce the
% dimensionality, a POD approach is adopted. The matrix $\Xm$ is decomposed using
% the singular value decomposition as:
\begin{equation}
\Xm = \Um \Sm \Vm^*,
\end{equation}
where the matrix $\Um \in \mathbb{R}^{N \times (m-1)}$ contains the orthogonal left-singular vectors. We can
then project the operator onto the space spanned by the left-singular vectors
to get the reduced operator $\ATm$. It is possible to note that the reduced
operator does not require the construction of the high-dimensional one:
\begin{equation}
\ATm = \Um^*\Am\Um = 
    \Um^*\Ym\Xm^\dagger \Um = 
    \Um^* \Ym \Vm \Sm^{-1} \Um^* \Um = 
    \Um^* \Ym\Vm \Sm^{-1}.
\end{equation}
We can now reconstruct the eigenvectors and eigenvalues of the matrix
$\mathbf{A}$ thanks to the eigendecomposition of $\mathbf{\tilde{A}}$ as
$\mathbf{\tilde{A}} \mathbf{W} = \mathbf{W} \bm{\Lambda}$. In
particular each nonzero eigenvalue $\lambda$ in $\bm{\Lambda}$ is a
DMD eigenvalue. The corresponding DMD eigenvectors, the so called
\textit{exact} modes~\cite{tu2014dynamic}, can be retrieved by
the eigenvectors of $\mathbf{\tilde{A}}$ as $\PHIm = 
\mathbf{Y}\mathbf{V} \Sm^{-1} \mathbf{W}$, where different scalings are
possible. We underline that each pair $(\phi, \lambda)$ computed as
above is an eigenpair of $\mathbf{A}$ (please refer to the proof of
Theorem 1 in~\cite{tu2014dynamic}). Thus, being $\Am = \PHIm
\LAMm \PHIm^\dagger$, we can approximate the evolution of the system $\xv_{k+1}
= \PHIm\LAMm\PHIm^\dagger\xv_k$. Moreover, it is easy to demonstrate that the
approximation of a generic future snapshots can be computed as:
\begin{equation}
    \xv_{k+j} = \PHIm \LAMm^j \PHIm^\dagger\xv_k.
\end{equation}

In this work we compute the DMD modes of the matrix composed by
the value of the time-varying lift coefficient for a set of given
geometrical parameters. Then we can predict the future state of the
coefficient and, using a regression method, approximate the target
function at untried new parameters. All the DMD computation have been
carried out by the Python package PyDMD~\cite{demo18pydmd}.

\section{Global sensitivity analysis through Active Subspaces}
\label{sec:as}
Active subspaces~\cite{constantine2015active} have been successfully
employed in many engineering
fields~\cite{constantine2014active,constantine2015exploiting}. Among
other we mention applications in shape
optimization~\cite{lukaczyk2014active,ghoreishi2019adaptive},
combustion simulations~\cite{ji2018shared}, in naval
engineering~\cite{tezzele2018dimension}, and in optimization procedures
coupled with genetic algorithm~\cite{demo2020asga}.
For multifidelity dimension reduction with AS
see~\cite{lam2018multifidelity}, for multivariate extension of AS we
mention~\cite{zahm2018gradient}, while for a coupling with deep neural
networks see~\cite{tripathy2019deep}. A kernel-based extension of AS
for multivariate functions we be found in~\cite{romor2020kas}.

Active subspaces have also been proven as a useful tool to enhance
model order reduction techniques such as proper orthogonal
decomposition (POD) with interpolation for structural and fluid dynamics
problems~\cite{demo2019cras}, and POD-Galerkin methods for a
parametric study of carotid artery stenosis~\cite{tezzele2018combined}.

Here we briefly introduce the active subspaces property for functions
not depending on time, for the
details and estimates regarding the method we refer
to~\cite{constantine2015active}. For the actual computations to find
AS we used the open source Python package ATHENA - Advanced
  Techniques for High dimensional parameter spaces to Enhance
  Numerical Analysis~\cite{athena}, derived in part from the Python Active
subspaces Utility Library~\cite{constantine2016python}. 

Let $\mupar \in \mathbb{R}^k$ the parameters of our problem, $f$ be a
parametric scalar function of interest $f(\mupar): \mathbb{R}^k \to
\mathbb{R}$, and $\rho: \mathbb{R}^k \to \mathbb{R}^+$ a probability
density function representing uncertainty in the input parameters. Active subspaces
are a property of the pair $(f, \rho)$. They are defined as the leading
eigenspaces of the second moment matrix of the target function's
gradient and constitutes a global sensitivity index more general than
coordinate-aligned derivative-based ones~\cite{zahm2018gradient}.

The second moment matrix of the gradients $\mathbf{C}$, also called uncentered 
covariance matrix of the gradients of $f$ with
respect to the input parameters, is defined~as
\begin{equation}
\mathbf{C} = \mathbb{E}\, [\nabla_{\mupar} f \, \nabla_{\mupar} f
^T] =\int (\nabla_{\mupar} f) ( \nabla_{\mupar} f )^T
\rho \, d \mupar,
\end{equation}
where $\mathbb{E} [\cdot]$ is the expected value, and $\nabla_{\mupar} f
\equiv \nabla f({\mupar}) \in \mathbb{R}^k$. $\mathbf{C}$ is symmetric
thus it admits a real eigenvalue decomposition that reads: 
\begin{equation}
\mathbf{C} = \mathbf{W} \Lambda \mathbf{W}^T,
\end{equation}
where $\mathbf{W}$ indicates the orthogonal matrix containing the eigenvectors
of $\mathbf{C}$ as columns, and $\Lambda$ is a diagonal matrix
composed by the non-negative eigenvalues arranged in descending order.
We can decompose the two matrices as follows
\begin{equation}
\Lambda =   \begin{bmatrix} \Lambda_1 & \\
& \Lambda_2\end{bmatrix},
\qquad
\mathbf{W} = \left [ \mathbf{W}_1 \quad \mathbf{W}_2 \right ],
\qquad
\mathbf{W}_1 \in \mathbb{R}^{k\times M},
\end{equation}
where $M < k$ has to be properly selected by identifying a spectral
gap. In particular, we define the active subspace of dimension $M$ as
the principal eigenspace corresponding to the eigenvalues prior to the
gap. Then we can map the full parameters to the reduced ones through
$\mathbf{W}_1$. We define the active variable as $\mupar_M =
\mathbf{W}_1^T\mupar \in \mathbb{R}^M$, and the inactive variable as $\etapar =
\mathbf{W}_2^T\mupar \in \mathbb{R}^{k-M}$. In practice the matrix
$\mathbf{C}$ is constructed with a Monte Carlo procedure.

AS stipulates that the directional derivatives in directions belonging
to the kernel of $\mathbf{W}_1^T$ are significantly smaller that those
belonging to the range of the same matrix. Moreover this assumptions
are made in expectation rather then in absolute
sense~\cite{wycoff2019sequential}.

Since in this way we are considering a linear combinations of the input
parameters, we can associate the eigenvectors elements to the weights
of such combinations, thus providing a sensitivity of each
parameter. We underline that if a weight is almost zero, 
that means $f$ does not vary along that direction on average.

We can use the active variable to build a ridge function
$g$~\cite{lin1993fundamentality} to approximate
the function of interest, that is
\begin{equation}
f (\mupar) \approx g(\mathbf{W}_1^T\mupar) = g (\mupar_M).
\end{equation}

In this work we want to study the behaviour of a target function $f(\mupar, t):
\mathbb{R}^k \times \mathbb{R}^+ \to \mathbb{R}$
that depends on the parameters $\mupar$ and on time $t$ as well. This
results in extending the active subspaces property to dynamical
systems, that means having to deal with time-dependent uncentered
covariance matrix $\mathbf{C} (t)$, and corresponding eigenvectors
$w_i(t)$. Efforts in this direction has been done
in~\cite{constantine2017time} for a lithium ion battery model,
in~\cite{loudon2016mathematical} for long term model of HIV infection
dynamics, and more recently an application of dynamic mode
decomposition and sparse identification to approximate one-dimensional active
subspaces in~\cite{aguiar2018dynamic}. In these works they refer to
dynamic active subspaces (DyAS) as the time evolution of the active
subspaces of a time-dependent quantity of interest.

DyAS are useful to assess the importance of each input parameter at
given times and to study how the weights associated to the inputs
evolve. In the following we are going to compute the AS for a set of
equispaced times $t_i$. If some of the parameters are almost zero in
the entire time window we can safely ignore them in the construction
of the Gaussian process regression.

\section{Computational pipeline}
\label{sec:results}

\begin{figure}[t]
\centering\includegraphics[width=.95\textwidth]{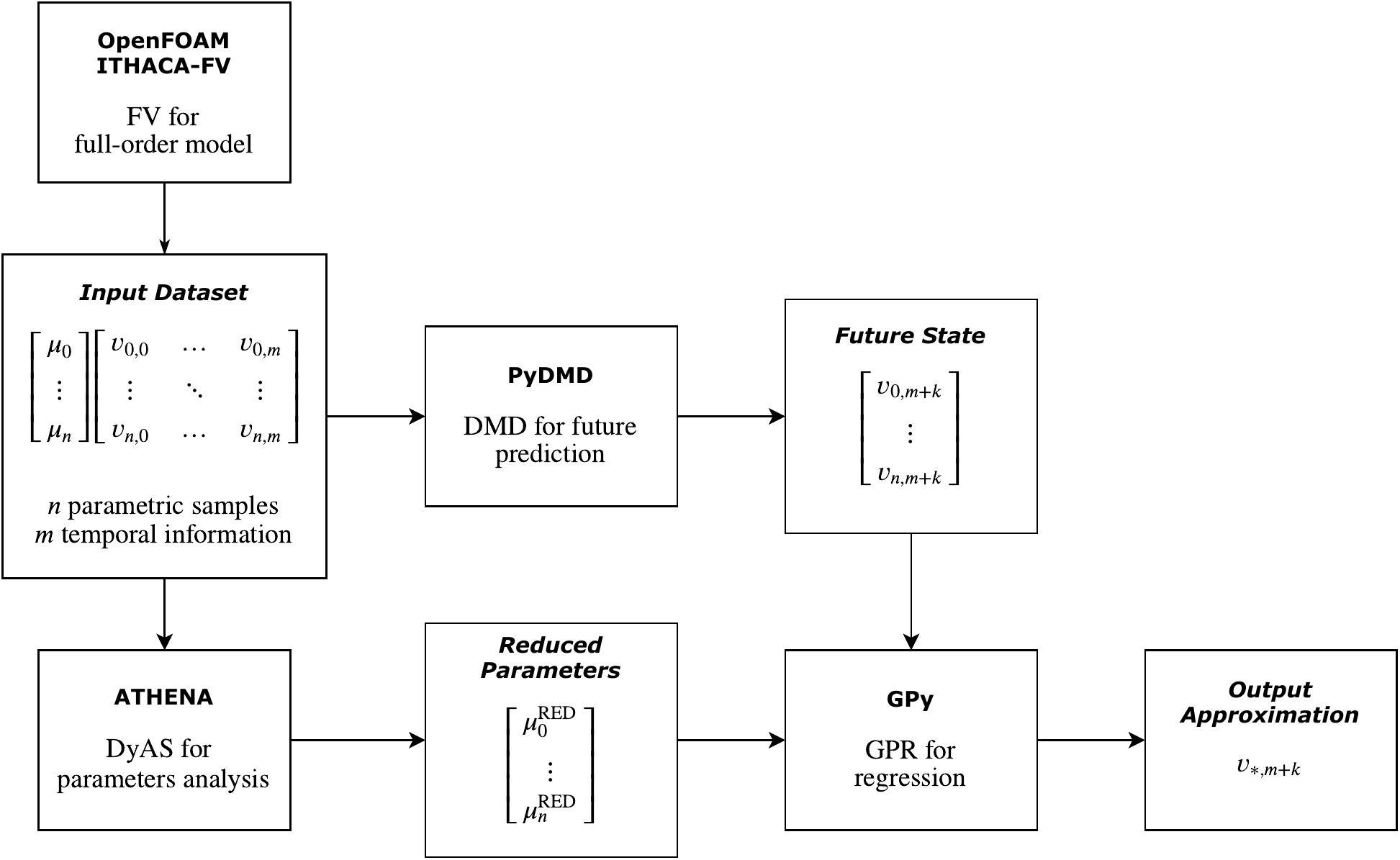}
\caption{Flowchart representing the proposed computational
	pipeline.}\label{fig:pipeline}
\end{figure}

In the present section we will discuss the numerical experiments carried out
to test the DyAS analysis and present the results obtained. As
reported in \autoref{sec:problem}, each high fidelity simulation is based on a
parametric fluid dynamic model governed by the Reynolds Averaged
Navier--Stokes (RANS) equations. Thus, a number of flow simulations
have been carried out selecting different samples in the parametric space to
test the performance --- in terms of lift coefficient --- of different airfoil shapes.
The simulations made use of both the RANS solver provided in the OpenFOAM~\cite{of}
finite volumes library, and of the DMD acceleration methodology
described in \autoref{sec:dmd}. Once the lift coefficients output were available for
all the samples tested in the input parameters space, the  DyAS analysis was applied
to assess possible parameter redundancy. The elimination of the redundant parameters
detected in the DyAS analysis allowed for the generation of a surface response model
based on a lower dimensional space, which has been finally tested against the original
RANS model accelerated through DMD, and against the surface response model
based on the original input parameter
space. \autoref{fig:pipeline} graphically summarizes the proposed
  pipeline, clarifying how the methods (and the software) are
  integrated together, while the following sections will further detail 
each part of the computational pipeline just outlined.

\subsection{Parametric shape deformation}
\begin{figure}
\begin{minipage}{0.53\textwidth}
\centerline{
  \ifpdf
  \resizebox{1\textwidth}{!}{
    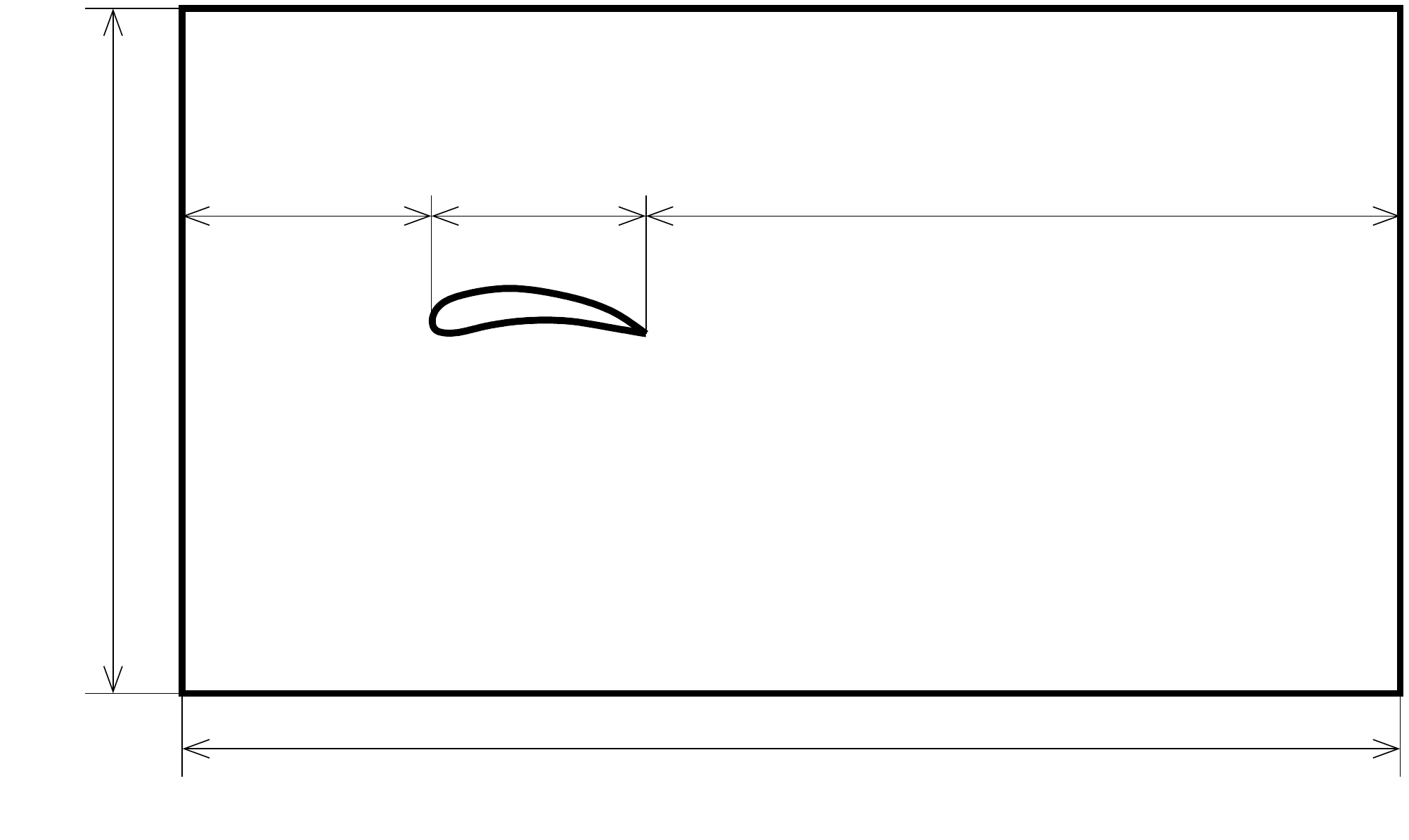
  }
  \else
  \resizebox{1\textwidth}{!}{
    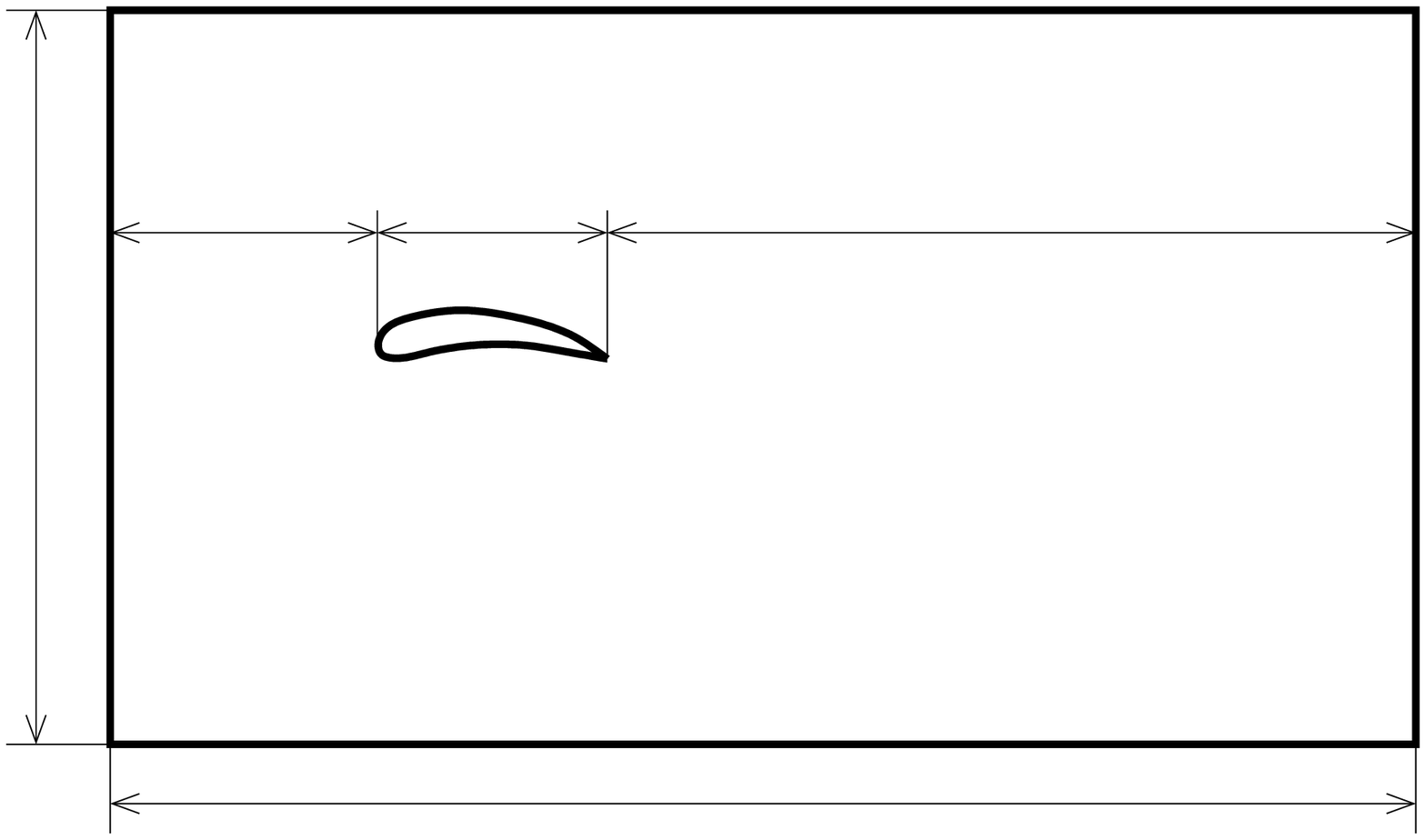
  }
  \fi
}
\end{minipage}
\begin{minipage}{0.47\textwidth}
\includegraphics[width=\textwidth]{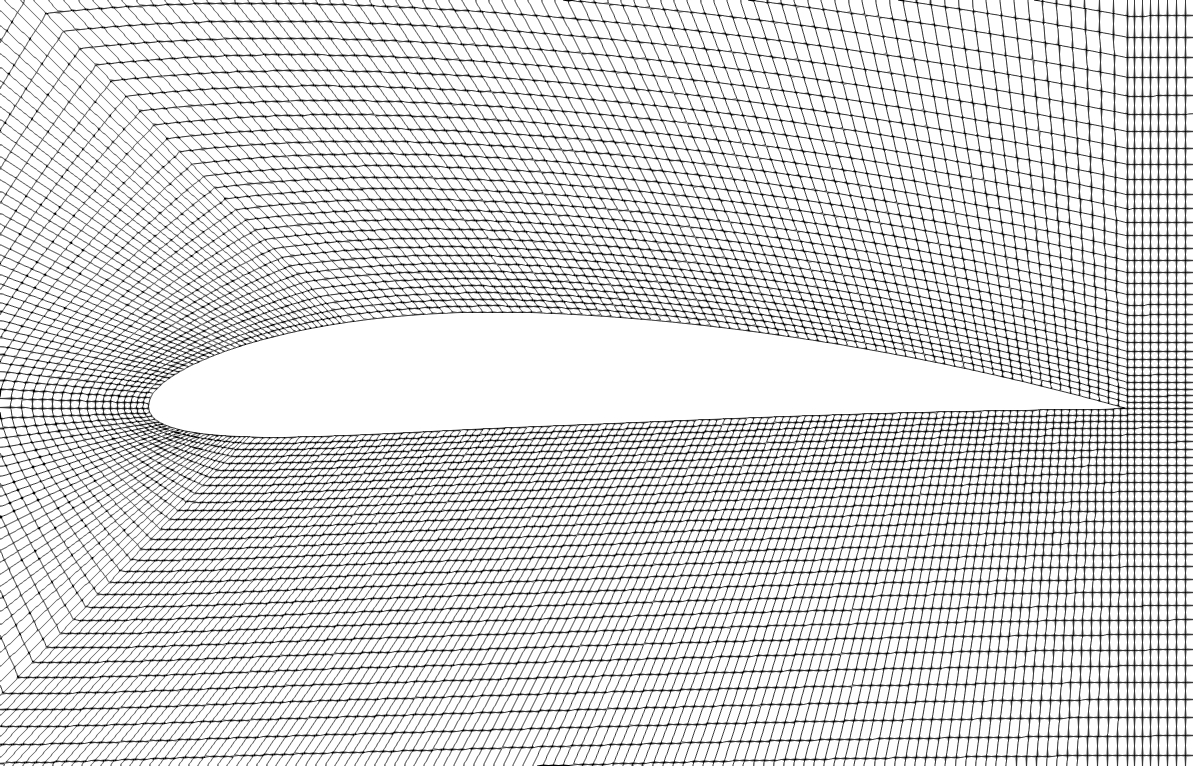}
\end{minipage}
\caption{Sketch of the computational domain used to solve the fluid
  dynamics problem in its reference configuration. The left picture
  reports a schematic view on the domain with the main geometrical
  dimensions. The right plot reports a zoom on the mesh in the
  proximity of the wing.}
\label{fig:mesh}
\end{figure}

\begin{figure}
\centering
\begin{minipage}{0.49\textwidth}
\includegraphics[width=\textwidth]{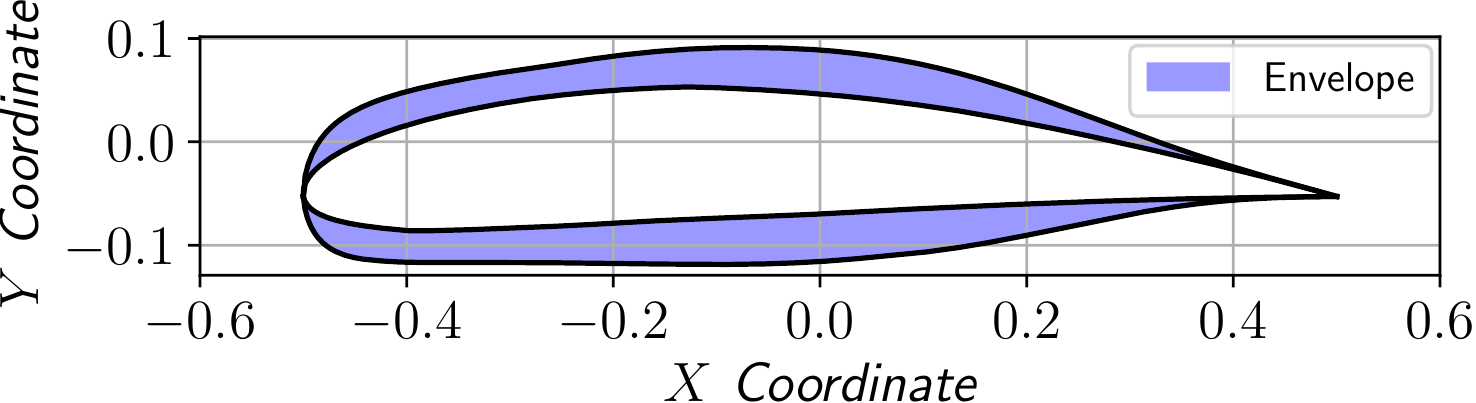}
\end{minipage}
\begin{minipage}{0.49\textwidth}
\includegraphics[width=\textwidth]{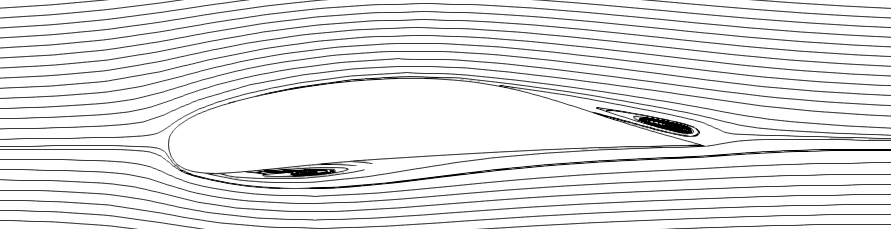}
\end{minipage}
\caption{The left picture reports in light blue the envelope of all
  the tested configurations used during the training stage. The right
  picture depicts the flow velocity streamlines for one particular
  sample inside the training set $\bm \mu = [0.0071; 0.0229; 0.0015;
  0.0015; 0.0087; 0.0107; 0.0033; 0.0130; 0.0247; 0.0280]$.}
\label{fig:result}
\end{figure}

The fluid dynamics problem is resolved using the finite volume method.
The wing is immersed in a rectangular domain according to \autoref{fig:mesh}.
The reference mesh counts $46500$ hexahedral cells
and is constructed using the \emph{blockMesh} utility of the
OpenFOAM library. \autoref{fig:mesh} depicts a detail of the grid
in proximity of the wing. The meshes in the deformed configuration
have been obtained starting from the reference configuration using a
radial basis function smoothing algorithm similar to the one
implemented in~\cite{Bos2013}. A single deformation corresponds to a
sample $\mupar$ in the parameter space  $\mathbb{D} := [0, 0.03]^{10} \subset
\mathbb{R}^{10}$.  Therefore all the deformed meshes share
the same number of cells and the same mesh topology.  In particular
Wendland~\cite{Wendland1995} second order kernel functions with radius
$r_{\text{RBF}} = 0.1~\si{m}$ have been used. The control points of the RBF
procedure have been placed on each mesh boundary point located onto the wing
surface. Since the outer boundary points are fixed we decided to
neglect them from the RBF computation using a smoothing function
defined in such a way that the RBF contribution reduces to zero after
a certain distance from a focal
point~\cite{Jakobsson2007}. Particularly, the focal point has been
placed in the geometric center of the airfoil chord segment and the distance
from the focal point after which the RBF contribution is neglected is set
to $r_{\text{out}} = 7~\si{m}$. In \autoref{fig:result} we depict the
envelope of all the tested configurations, and the flow velocity
streamlines for a particular sample in the parameter space. 
A uniform and constant velocity equal to $\bm
u_{\text{in}}=1~\si{m}/\si{s}$ is set at the inlet boundary, while the
constant value of the kinematic viscosity is set to $\nu = 2\mathrm{e}{-5}~\si{m^2}/\si{s}$.
This configuration, considering a chord length $D=1~\si{m}$, corresponds to Reynolds number
$\mathrm{Re}=50000$. As well known, a flow characterized by Reynolds number of such
magnitude requires turbulence modeling to be numerically simulated with reasonable
computational effort. In the present work, turbulence has been modeled using a RANS
approach with a Spalart-Allmaras turbulence
model~\cite{SPALART1992}. The pressure velocity coupling is resolved
in a segregated manner making use of the PIMPLE algorithm which merges the
PISO~\cite{Issa1986} and the SIMPLE~\cite{Patankar1972} algorithm. The
time step used to advance the simulation in time is set constant and
equal to $\Delta t=1\mathrm{e}-3~\si{s}$. The convective terms have
been discretized using a second-order upwinding scheme, while the
diffusion terms are discretized using a linear approximation scheme
with non-orthogonal correction. The time discretization is resolved
using a second order backward differentiation formula. The simulation
is advanced in time until the flow has reached stationary behavior.
For the present problem, setting a total simulation time $T_s = 30~\si{s}$ is
sufficient to reach a solution which is reasonably close to the steady state one.
In order to check the consistency of the numerical results, the stationary lift 
coefficient computed for the reference configuration, which corresponds to a 
standard NACA $4412$ profile with a $0^{\circ}$ angle of attack, has been 
compared with data from literature \cite{airfoiltools}. The computed lift 
coefficient for such setting is equal to $C_L = 0.355$ and the available 
reference value varies between $C_L=0.1804$ and $C_L=0.3708$ depending on the 
value of $N_{\mbox{crit}}$ (which is used to model the turbulence of
the fluid or roughness of the airfoil). Therefore, our numerical results are in line with 
available data in existing literature\footnote{Such comparison is not 
exhaustive to completely verify the accuracy and the reliability of the full 
order model numerical simulations. It is however beyond the scope of this work 
to perfectly match experimental activities or previous numerical results with 
the full order simulations. More accurate FOM results would of course result in 
more accurate ROM results but would not affect the presented methodology.}.

\subsection{Parameter space reduction}
The present section will discuss the application of DyAS to the problem of the
two dimensional turbulent flow simulation past airfoil sections with parameterized shape.
Such a fluid dynamic problem is relevant in several engineering
fields, as it is encountered in a number of industrial applications,
ranging from aircraft and automotive design, to turbo machinery and
propeller modeling. We must here point out that in this work, the DMD
method is used for faster evaluation of the parameterized airfoils
lift towards a steady state regime solution. We remark that, since DMD is designed for
time evolutionary problems, the same procedure can be used in the same fashion,
to speed up convergence to periodic regime solutions~\cite{le2015flow}.
Indeed, recent work on hydroacoustic computations based on LES
suggested that DMD modal decomposition can successfully be employed in
the reconstruction of complex and turbulent flow fields
\cite{Gadalla2020}  provided that the snapshots used are enough to
characterize all the relevant time and space frequencies in the
flow. In addition, we observed that complex full order flows
characterized by richer spectra require a higher amount of modes to
obtain accurate flow fields reconstruction. Thus, our experience
suggests that the ROM instruments used in this work are indeed
effective when employed with more complex physics. For such reason,
given our experience, we infer that the design pipeline here presented
can also be used to study the unsteady dynamics of bubbles and
vortices past the airfoil. Obviously one requirement of such type of
problems would be a suitable FOM model able to capture transition
phenomena occurring in the stall region. For example, we believe that
the underlying high fidelity URANS solver would not be appropriate and
that a transition to a LES approach would be required. For
projection-based ROMs in a turbulent setting
see~\cite{HiStaMoRo2019,StabileBallarinZuccarinoRozza2019,HijaziAliStabileBallarinRozza2020}.

\begin{figure}[h!]
\centering
\includegraphics[width=.9\textwidth, trim=70 0 70 0, clip]{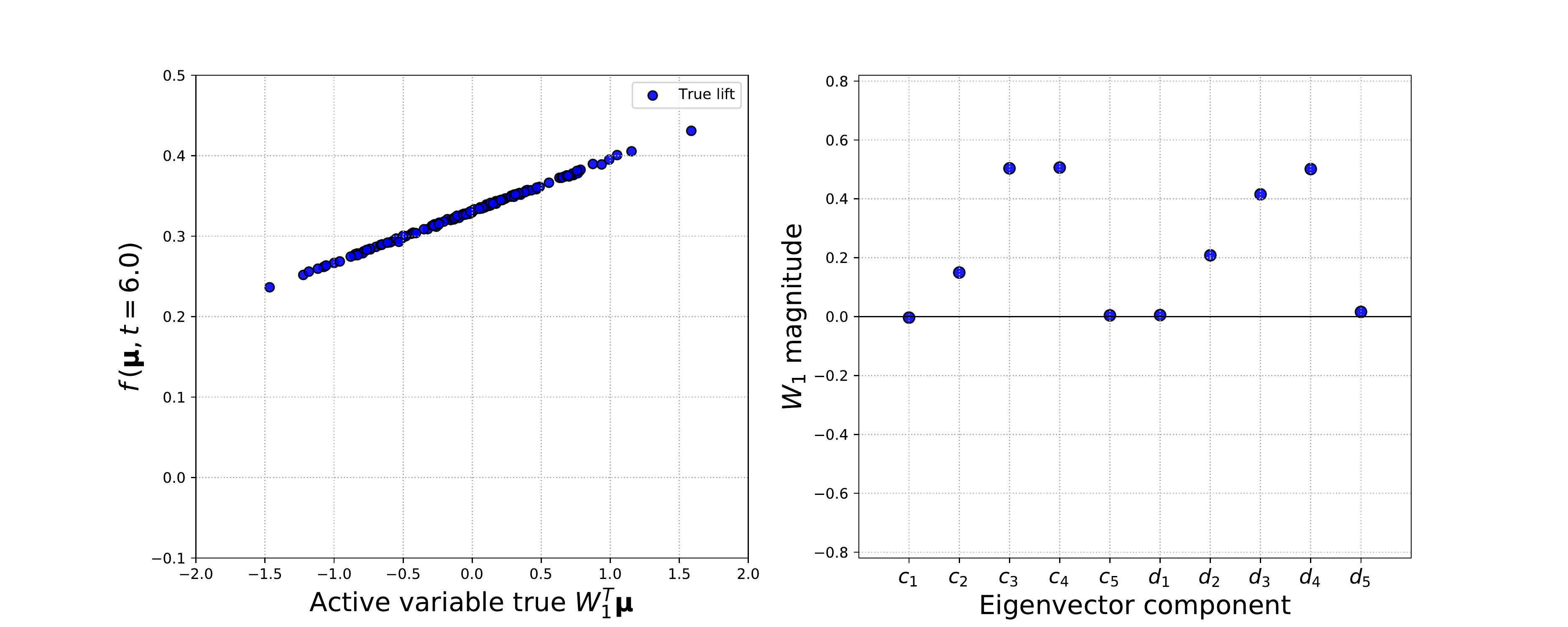}
\caption{On the left the sufficiency summary plot for the lift
coefficient at time $t=6.0$ seconds. On the right the first eigenvector
components at the corresponding parameters.}
\label{fig:ssp_lift_6s}
\end{figure}

\begin{figure}[h!]
\centering
\includegraphics[width=.9\textwidth, trim=70 0 70 0, clip]{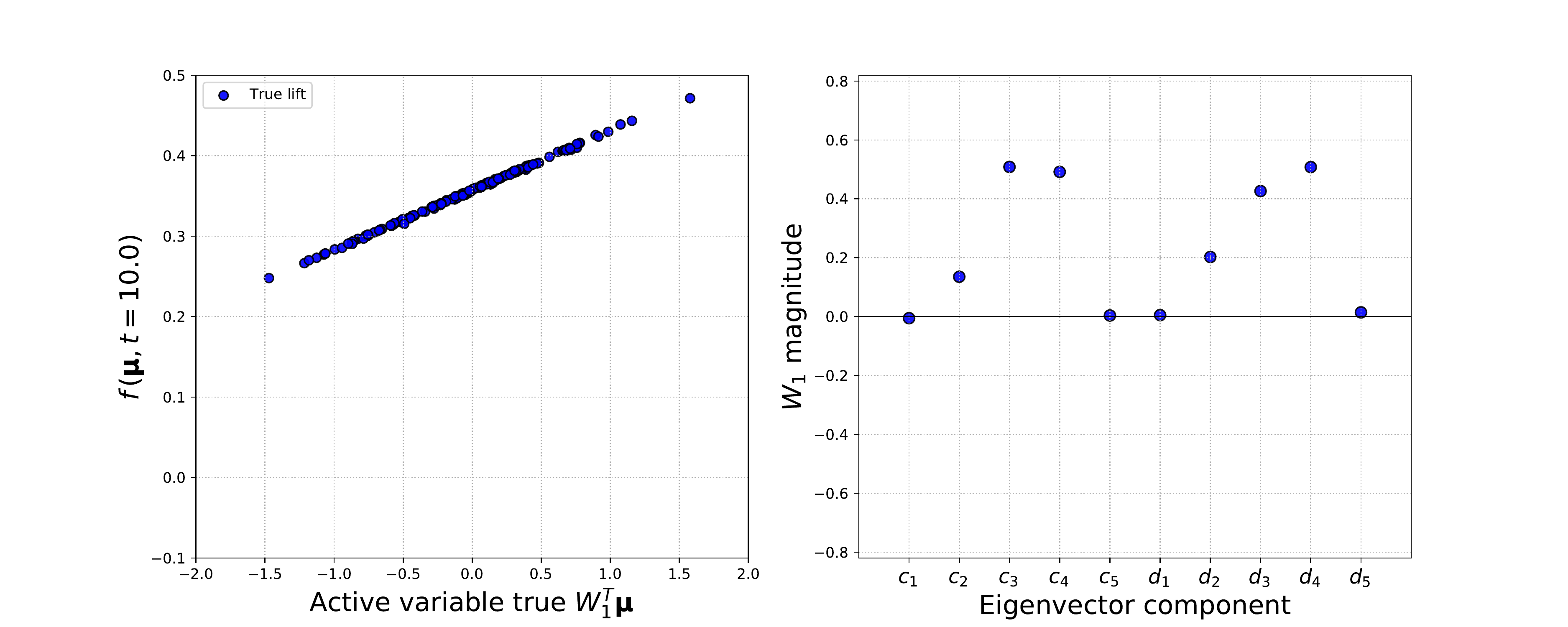}
\caption{On the left the sufficiency summary plot for the lift
coefficient at time $t=10.0$ seconds. On the right the first eigenvector
components at the corresponding parameters.}
\label{fig:ssp_lift_10s}
\end{figure}

\begin{figure}[h!]
\centering
\includegraphics[width=.9\textwidth, trim=70 0 70 0, clip]{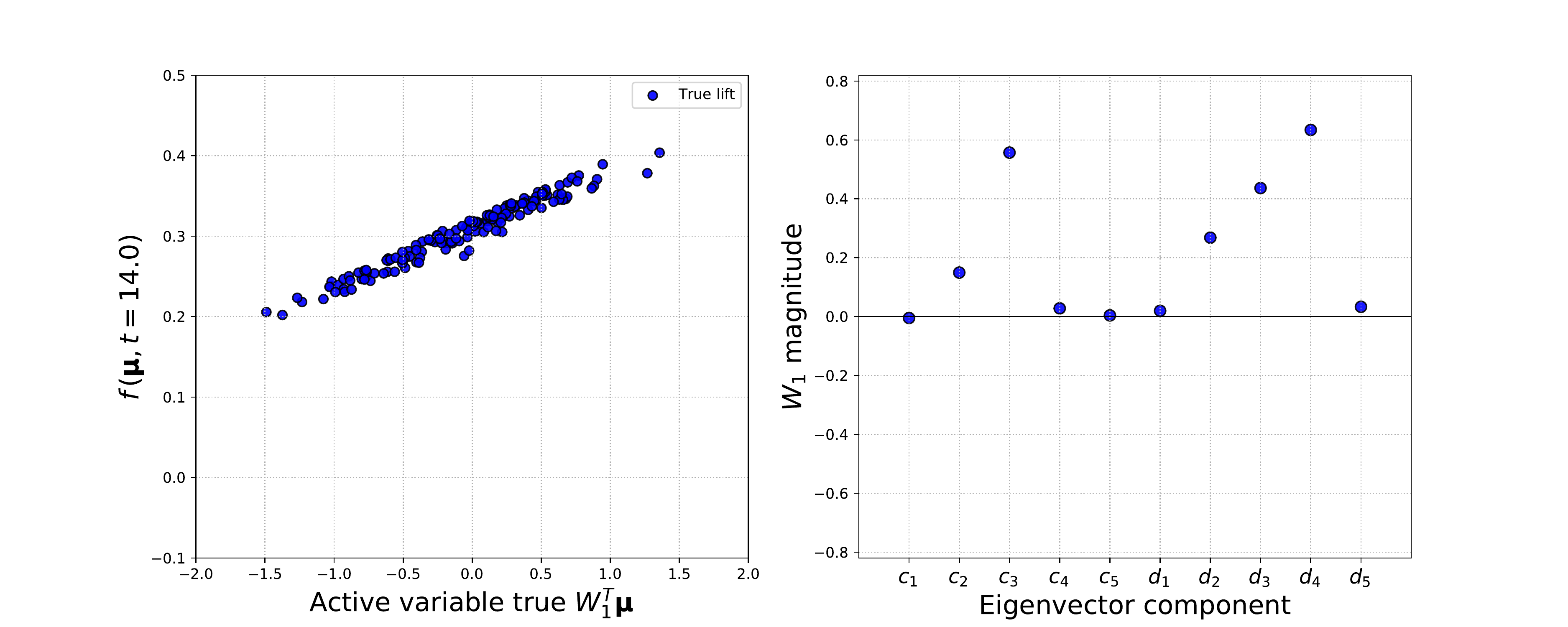}
\caption{On the left the sufficiency summary plot for the lift
coefficient at time $t=14.0$ seconds. On the right the first eigenvector
components at the corresponding parameters.}
\label{fig:ssp_lift_14s}
\end{figure}

\begin{figure}[h!]
\centering
\includegraphics[width=.9\textwidth, trim=70 0 70 0, clip]{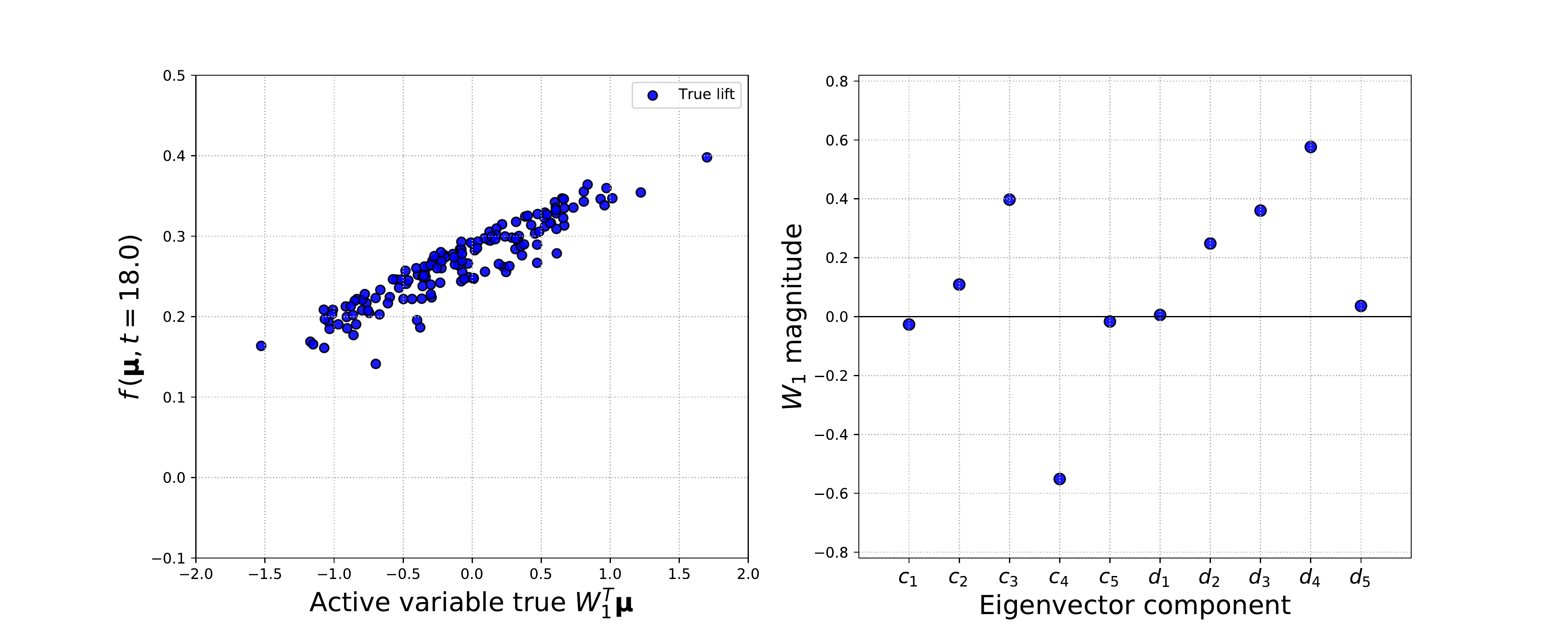}
\caption{On the left the sufficiency summary plot for the lift
coefficient at time $t=18.0$ seconds. On the right the first eigenvector
components at the corresponding parameters.}
\label{fig:ssp_lift_18s}
\end{figure}

A few plots describing the DyAS results for the lift coefficient
output are presented
in~\autoref{fig:ssp_lift_6s},~\ref{fig:ssp_lift_10s},~\ref{fig:ssp_lift_14s},
and~\ref{fig:ssp_lift_18s}. The plots in the figures are aimed at
representing the evolution of the active subspace effectiveness and
composition over the time dependent flow simulations. More
specifically, the left diagram in each figure plots the lift
coefficient at each sample point tested, as a function of the first
active variable obtained through a linear combination of the sample
point coordinates in the parameter space, that is $f(\mupar, t)$
against $\mathbf{W}_1^T \mupar$.
\begin{figure}[h!]
\centering
\includegraphics[width=.7\textwidth]{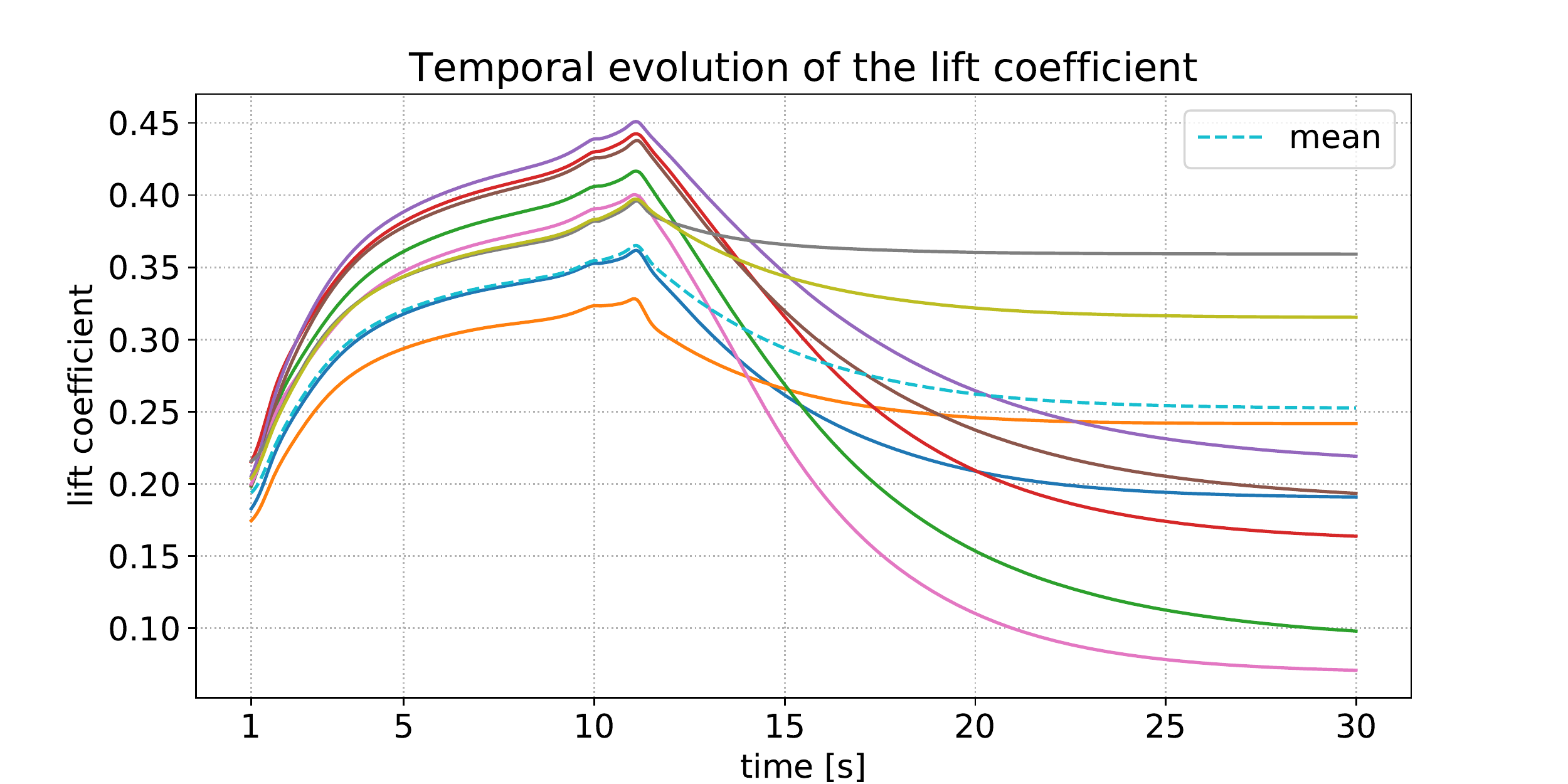}
\caption{The temporal evolution of the lift coefficient from
  $1~\si{s}$ to $30~\si{s}$ for $9$
different parameters, together with the mean (dashed). The angle
of attack is fixed for all the airfoil profiles and it is equal to $0^\circ$.}
\label{fig:temporal_lift}
\end{figure}
Presenting the components of the
first eigenvector of the uncentered covariance matrix, the right plot
in each figure indicates the weights used in such linear combination
to obtain the first active variable. In summary, the right diagram in
each Figure suggests the impact of each of the original parameters on
the first active variable, while the left diagram is an indicator of
how well a one dimensional active subspace is able to represent  the
input to output relationship. Following the evolution of these two
indicators it is possible, at each time instant, to assess how
effective the one dimensional parameter dimension reduction is, and
what is the sensitivity of the reduced lift coefficient output to
variations of the original parameters.
The plots in ~\autoref{fig:ssp_lift_6s},~\ref{fig:ssp_lift_10s},~\ref{fig:ssp_lift_14s},
and~\ref{fig:ssp_lift_18s} show the results of the DyAS at the fixed time
instants $t=6~\si{s}, 10~\si{s}, 14~\si{s}, 18~\si{s}$, respectively.
We here remark that, given the aforementioned considerations about the solution
build up in the first 12 seconds of the simulations, the solutions at $t=6~\si{s}$ and $t=10~\si{s}$ 
are not entirely relevant by a physical perspective. Yet, presenting such cases is still helpful
in illustrating how the DyAS evolve over time and can be used to
evaluate the system behavior and the output sensitivities with respect
to the input parameters. For completeness in
\autoref{fig:temporal_lift} we depicted the temporal evolution of $9$
different morphed airfoils, and the mean among all the airfoils.
A first look at the right plots
for each time steps, suggests that the contribution of the parameters
corresponding to the bump shape functions
$r_1$, and $r_5$, for both the top and the bottom part of the airfoil
profile are almost negligible. This means the lift coefficient is
almost insensitive to variations of these 4 parameters. Alternatively,
it can be said that the output function is on average almost flat along directions corresponding
to the axes corresponding to parameters $c_1$, $c_5$, $d_1$, and $d_5$.

\autoref{fig:ssp_lift_6s} and~\ref{fig:ssp_lift_10s} present the
characterization of the one dimensional active subspace at time
$t=6~\si{s}$ and $t=10~\si{s}$, respectively. We can clearly see that
the lift coefficient is perfectly approximated along the identified direction,
and such direction (the eigenvector elements) is almost the
same at $t=6~\si{s}$ and $t=10~\si{s}$. This should not completely surprise as
both time instants are included in an initial acceleration phase
during which the air coming from the inflow boundary is reaching the
airfoil. Given the domain arrangement described in~\autoref{fig:mesh},
the flow velocity around the impulsively started airfoil leading edge  
is expected to reach the inflow value at time $t=10~\si{s}$. For such
reason, we will focus the description on the plots for $t=10~\si{s}$,
although the considerations can be immediately reproduced for previous
time steps. The left plot in~\autoref{fig:ssp_lift_10s} suggests that at
this meaningful instant, the first active subspace represents the
input to output relationship with remarkably good accuracy. In fact,
only a single output value corresponds to each active variable
value. In other words, when plotted against the first variable, the
output appears like a curve --- a line in the present case.  A look at
the right diagram suggests that the shape parameters having the most
impact on the lift generated by the airfoil are $c_3$, $c_4$, $d_3$
and $d_4$, which are the ones associated to shape functions with peaks
located around the middle of the airfoil chord. The positive values of
the eigenvector components associated to $c_3$, $c_4$, $d_3$ and
$d_4$, along with the positive slope of the curve in the left plot
in~\autoref{fig:ssp_lift_10s} suggest that, at this particular
time instant, higher values of lift can be obtained by increasing the  
airfoil thickness in the mid-chord region. 

Similar considerations can be drawn from \autoref{fig:ssp_lift_14s}, which
refers the the DyAS analysis carried out at $t=14~\si{s}$. Here, the points in
the left diagram do not completely cluster on top of a single valued curve
as was the case for the previous time step considered. Compared to what
has been observed at  $t=10~\si{s}$, the data clearly indicate that
at $t=14~\si{s}$ an input to output relationship obtained using only a one
dimensional active subspace will lead to less accurate lift coefficient
predictions.
Yet, the points in the plot are still all located within a
rather narrow band surrounding a regression line having positive slope.
Thus, all the considerations on the lift coefficient sensitivity with
respect to variations of the shape parameters that can be inferred from
the right plot, will still hold at least by a qualitative standpoint. Here,
the eigenvector components suggest that the most influential parameters on
the lift coefficient are  $c_3$, $d_3$ and $d_4$, while $c_2$
and $d_2$ affect the output in lesser but not negligible
fashion. Compared to the previous case the importance of coefficient $c_4$ 
on the output is significantly reduced. We recall that $c_4$ is associated with
increased $y$ coordinates of the airfoil suction side past the mid-chord region.
Thus, we might infer that in the acceleration phase higher lift values are
obtained not only increasing the front thickness, but also lowering the camber
line in the region past mid-chord.

\autoref{fig:ssp_lift_18s} shows the results of the DyAS analysis at $t=18~\si{s}$,
when the flow approaches the final regime solution. Following the trend observed
for $t=14~\si{s}$, the left plot in the figure indicates that a one
dimensional active subspace is not
completely able to represent the input to output relationship in a satisfactory fashion.
With respect to the previous plots, the output values are here located in an even wider
band around a regression line with positive slope. Again, on one hand this increasingly
blurred picture suggests that higher dimensional active subspaces are required to reproduce
the steady state solution with sufficient accuracy; on the other
hand, the diagram still suggests a quite definite trend in the output, which can be
exploited for qualitative considerations. 
Quite interestingly, at the present time step the eigenvector
component corresponding to the $c_4$ coefficient has negative sign. Given the
positive slope of the input to output relationship in the left plot
of~\autoref{fig:ssp_lift_18s}, this implies that increases in the
airfoil ordinates on the top side in the region past the mid-chord
result in lift loss. Thus, this seems to suggest that an airfoil with
a higher camber line curvature, combined with a thicker leading edge
region might result in increased lift. This should not surprise, as a similar kind of
airfoil would result in a higher downwash due to the increased camber line curvature, yet
being able to avoid stall by means of a thicker and rounder leading edge.
Thus, the DyAS analysis at different time steps shows that as the impulsively started
airfoil moves from an acceleration phase to a steady state regime solution, the shape
modifications leading to increased lift transit from a purely symmetric increase of the
thickness in the mid-chord region, to a non-symmetric modification of the camber
line united with a symmetric leading edge thickness increase, respectively. Such
behavior is indicated by the sign of $c_4$ coefficient in the eigenvector characterizing the
one dimensional active subspace, which is likely detecting that at steady state, regime solution,
airfoils with higher camber line curvature and thicker leading edges
produced higher downwash.

We underline that the eigenvector components of all the time instants
presented corresponding to the coefficients $c_1$, $c_5$, $d_1$, and
$d_5$ are almost zero. This means that on average the lift coefficient
is almost flat along these directions. We are going to exploit this
fact by freezing these parameters and constructing a GPR on a reduced
parameter space.

\subsection{GPR approximation and prediction of the lift coefficient}

\begin{figure}
\centering
\includegraphics[width=1.\textwidth]{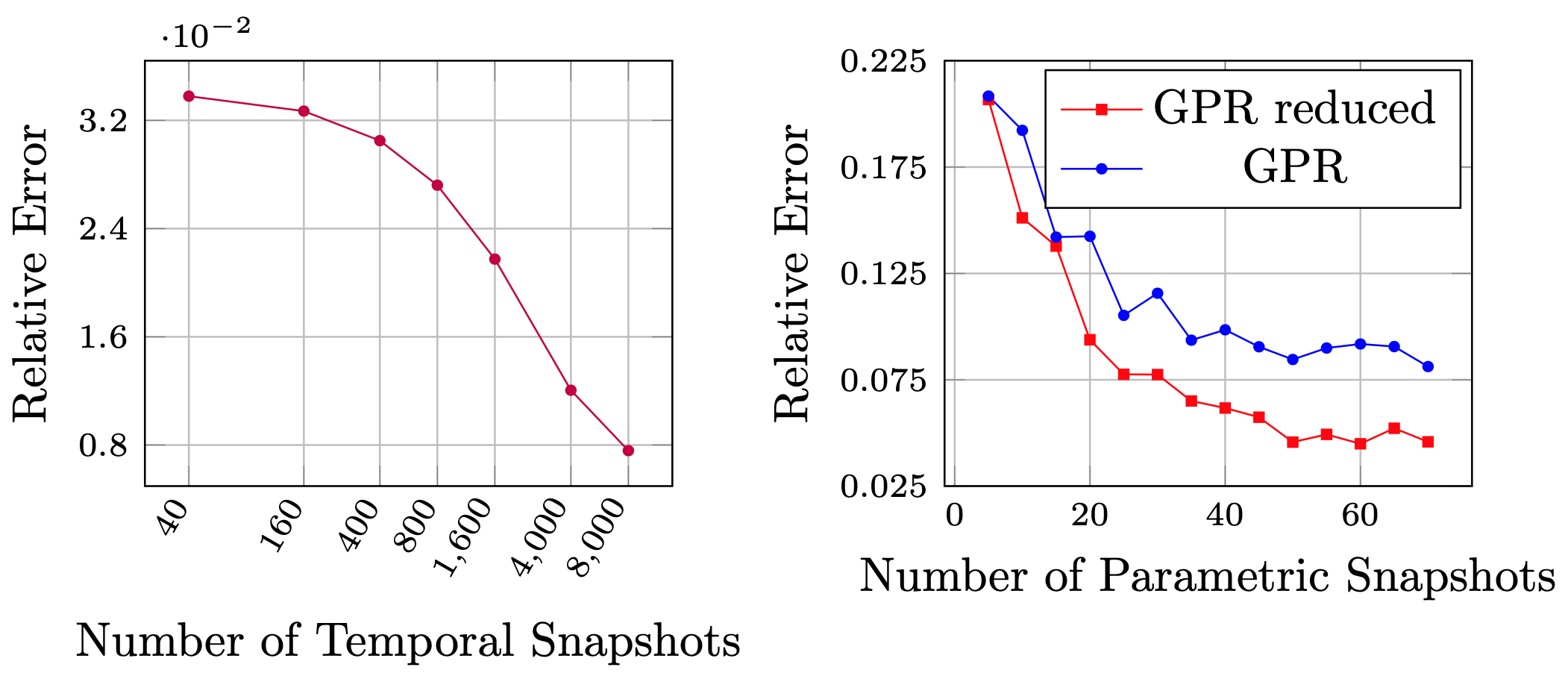}
    \caption{Sensitivity analysis of the dimension of the training set for
    the DMD (left) and for the response surface using GPR (right). For the
    DMD, we use 70 samples (of the parametric space) evolving in time in
    $[12, 20]~\si{s}$ and we measure the mean relative error at time
    $30~\si{s}$ varying the sampling frequency; for the GPR, we build the
    response surface using up to 70 sampling lift coefficients at time
    $20~\si{s}$ and computing the mean relative error over the test dataset
    composed by 100 test deformations.}\label{fig:sensibility}
\end{figure}

The previous analysis pointed out the presence of several input parameters
with minimal average influence on the target function. Making use of such consideration
we construct a response surface which only depends on the remaining parameters. Both
for the full parameter space and the reduced one, we use a Gaussian process regression with
a RBF kernel implemented in the open source Python package GPy~\cite{gpy2014}.
We then compare the performance of the two regression strategies by
computing the relative error over a test data set composed by 100
samples. The error is computed as the Euclidean norm of the difference
between the exact and the approximated solution over the norm of the
exact solution. The training set is composed by the same 70 samples,
in 10 dimensions for the GPR over the original parameter spaces, and
in 6 dimensions for the reduced one. Up to $t = 20~\si{s}$ the
training is done using the high-fidelity simulations.

\begin{figure}
\centering
\includegraphics[width=1.\textwidth]{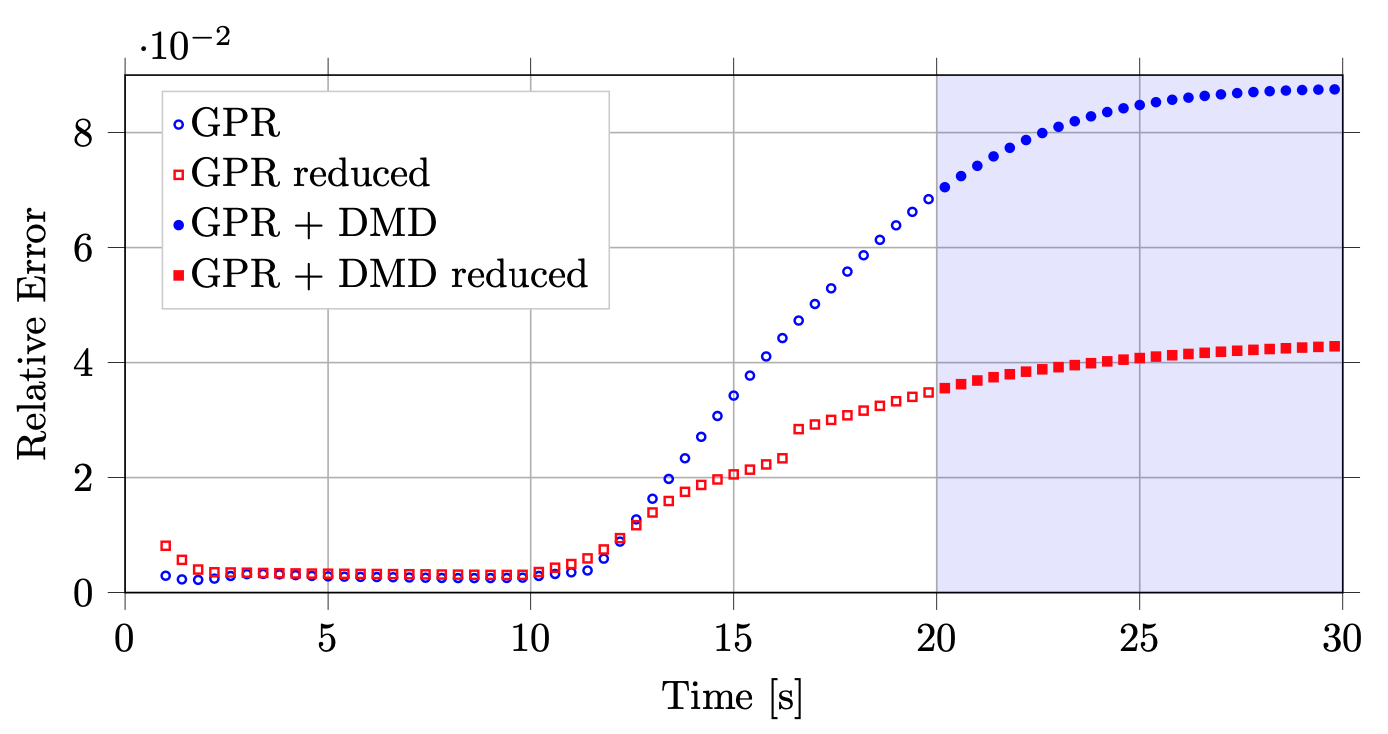}
\caption{The relative error of the approximated outputs at different times. The
relative error is computed on 100 test samples, using the
high-fidelity lift coefficient to
train the regression for $t \le 20~\si{s}$, while for $t > 20~\si{s}$ the DMD forecasted
states are used for the training.}\label{fig:error}
\end{figure}

To speed up the convergence to the regime state ($t = 30~\si{s}$) we applied
the DMD to get the future-state prediction of the lift. In particular, due
to the initial propagation of the boundary conditions, for all the 70 training
deformations we use the trend of lift coefficients within the temporal interval
$[12, 20]~\si{s}$ to fit the DMD model, that means $8000$ temporal information
($\Delta t = 0.001~\si{s}$). Since we used 10 POD modes --- selected
using the energetic criterion --- for the projection of the DMD operator, our
low-rank operator results of dimension $10$. Despite in this case the
dimensional reduction is not huge, this approach allows to predict the future
state in a very fast fashion. In the high-fidelity model, we need in fact
$1508$ CPU seconds (on average) to simulate $1$ second of the physical model,
instead using DMD we can approximate a future state in less than $0.1$ CPU
seconds. In practices, this means that, to reach the regime state with the
standard approach, the simulation lasts $1508~\si{s} \times 30 \approx
45000~\si{s}$, while with the DMD we have $1508~\si{s} \times 20 + 0.1~\si{s} \approx
30000~\si{s}$, guaranteeing to save $\textstyle\frac{1}{3}$ of the overall
computational load. All the simulations, both at the FOM and at the
ROM level have been run serially on an Intel Xeon E5-2640, 2.50GHz
CPU. We highlight that this is only a part of the computational saving
of the pipeline that we are proposing and is related to the training
stage. The DMD allows in fact for $1/3$ reduction of the simulation
time required to the FOM as the remaining time is simulated by an
approximated model. On the other side, once the reduced order model
has been constructed, exploiting the combination of the Gaussian
Process approximation and the DMD, it is possible to test new
geometries in real time, with a negligible computational cost. Regarding the accuracy, we present in
\autoref{fig:sensibility} a sensitivity analysis on the number of
training snapshots, varying the temporal sampling period $\Delta t_\text{DMD}$ from
$1\mathrm{e}-3~\si{s}$ to $0.2~\si{s}$ and measuring the error on the predicted
state at $t = 30~\si{s}$. Similarly, we propose an analysis on the GPR
accuracy: using a varying number of lift coefficients at $t = 20~\si{s}$, we build the
response surface and measure the error for untried parameters, both in the full dimensional space and in the reduced one.
In~\autoref{fig:error} we compare the two GPR performance at each of the time
steps analyzed in the simulations. Until $12~\si{s}$, the regressions behave in
a very similar fashion, while from $15~\si{s}$ the accuracy gain obtained by
distributing the 70 samples in a lower dimensional space becomes significant.
The error gap between the 6 and 10 dimensional response surface in fact,
consistently increases from $0.016$ at $15~\si{s}$ to $0.045$ at steady
state. This corresponds to a decrement of the error by a factor~$2$.

The proposed method achieves better results because it exploits the
DyAS to discard the directions of the input parameter space along
which the target function does not vary.

\section{Conclusions and perspectives}
\label{sec:the_end}

We presented a computational pipeline to improve the approximation of the
time-varying lift coefficient of a parametrized NACA airfoil. The
pipeline comprises automatic mesh deformation through RBF
interpolation, high-fidelity simulation with finite volume method of
turbulent flow past the airfoil, global sensitivity analysis
exploiting AS, and future state prediction via DMD reduced order
method. This resulted in more accurate Gaussian process regression
of the lift coefficient even if in a reduced parameter space.
Despite the turbulent nature of the flow, the selected testcase does not show
highly nonlinear phonemena --- e.g. stall, reattachment --- that usually
occur in several fluid dynamics problems. The proposed framework can be
extendend to address also more complex applications, provided that a suitable
number of snapshots is given to characterize the parameter space and
frequencies required by the DMD training. Of course such more demanding
training requirements would likely result in reduced ROMs speed up and would
require case-specific treatments.

After the creation of the high-fidelity solutions database the
application of AS highlighted a possible reduction of the parameter
space due to negligible contributions of $4$ different parameters. We
exploit this reduction to construct a GPR over a smaller parameter
space, thus improving its performance. Since the training of the
regression model is done over $6$ dimension instead of $10$, given the
same high-fidelity database dimension, the GPR is able to better
approximate the solution manifold. This results in better lift
coefficient predictions for new untried parameters. We also applied
DMD to have future-state prediction of the target function up to $30$
seconds and proved that the effective gain of the new GPR is preserved
also for any time after the $20$ seconds simulated with FV.
In particular from $13$ seconds the actual gain is significant, at $15$
seconds we have an increased performance by a factor~$2$ in the relative
error, which means that performing the regression in the reduced
parameter space produces a relative error equal to $0.02$, instead of
$0.036$. Evolving in the future the error drop increases up to
$0.045$ at regime ($0.042$ instead of $0.087$, keeping the factor~$2$).

This computational pipeline can be seen as a parametric dynamic mode
decomposition for some extent. Moreover, the sensitivity analysis has a
negligible computational cost with respect to the creation of the
offline high-fidelity database. 

Future developments can be the study of adaptive sampling
strategies exploiting a generic $n$-dimensional active subspace, and
the coupling of different model order reduction methods.
Another possible extension of the presented method regards the possibility
to apply the framework to a flow field --- e.g. pressure, velocity --- rather than
to a scalar output.
It would be interesting to
use this non-intrusive setting as a preprocessing tool to reduce the number of
simulations required to build a reduced basis space which is later used in an
intrusive manner~\cite{StaRo2018}.
We think this new computational pipeline can be of much interest in the
context of shape optimization and dynamical systems.

\section*{Acknowledgements}
This work was partially supported by an industrial Ph.D. grant sponsored
by Fincantieri S.p.A., and partially funded by the project UBE2 -
``Underwater blue efficiency 2'' funded by Regione Friuli Venezia Giulia,
POR-FESR 2014-2020, Piano Operativo Regionale Fondo Europeo per lo Sviluppo
Regionale.
It was also partially supported by European Union Funding for
Research and Innovation --- Horizon 2020 Program --- in the framework
of European Research Council Executive Agency: H2020 ERC CoG 2015
AROMA-CFD project 681447 ``Advanced Reduced Order Methods with
Applications in Computational Fluid Dynamics'' P.I. Gianluigi Rozza.

% \bibliographystyle{abbrv}
% \bibliography{mortech_biblio}

\end{document}